\documentclass{amsart}
\input amssym.def\input amssym

\newtheorem{thm}{Theorem}
 \newtheorem{cor}[thm]{Corollary}
 \newtheorem{lem}[thm]{Lemma}
 \newtheorem{prop}[thm]{Proposition}
 \newtheorem{defn}[thm]{Definition}

\newcommand{\ma}{\mathcal{A}}
\newcommand{\mb}{\mathcal{B}}

\newcommand{\mm}{\mathcal{M}}

\newcommand{\ha}{H_\mathcal{A}}

\begin{document}

\title{Rank-preserving module maps }

\author{Bin Meng}
\address{College of Science, Nanjing University of Aeronautics and Astronautics, Nanjing 210016,  People's Republic
 of China}

 \email{b.meng@nuaa.edu.cn}


\date{}


\begin{abstract}In this paper, we characterize rank one preserving
module maps on a Hilbert $C^\ast-$module and study its applications
on free probability theory.
\end{abstract}

\maketitle

\section{Introduction and
preliminaries}

The study of linear maps on operator algebras that preserve certain
properties has attracted the attention of many mathematicians in
recent decades. They have been devoted to the study of linear maps
preserving spectrum, rank, nilpotency , etc. In our study of free
probability theory, we find that module maps on Hilbert
$C^\ast-$module preserving certain properties are also important
(see \cite{meng1,meng2}) and thus the study of the module maps'
preserver problem becomes attractive.  A (left) Hilbert
$C^\ast-$module over a $C^\ast-$algebra $\ma$ is a left $\ma-$module
$\mm$ equipped with an $\ma-$valued product
$\langle\cdot,\cdot\rangle$ which is $\ma-$linear in the first and
$\ma-$conjugate linear in the second variable such that $\mm$ is a
Banach space with the norm $\|v\|=\|\langle
v,v\rangle\|^{\frac{1}{2}}, \forall v\in\mm$. Hilbert
$C^\ast-$modules  are introduced and first investigated in
\cite{kap} and a good textbook is \cite{lan}.

Rank-preserving problem is a basic problem in the study of linear
preserver problem (see \cite{hou2},\cite{hou1}). So we believe that
it should be the key to the study of modular preserver problems.
Rank-Preserving linear maps have been studied intensively by Hou in
\cite{hou1}.

Hilbert $C^\ast-$modules first appeared in the work of Kaplansky
\cite{kap}, who used them to prove that derivations of type I
$AW^\ast-$algebras are inner. He generalized Hilbert space inner
product to the values in a commutative unital $C^\ast-$ algebra. Let
$H$ be a separable infinite dimentional Hilbert space and let $\ma$
be a commutative unital $C^\ast-$ algebra. The Hilbert $\ma-$ module
$H\otimes \ma$ plays a special role in the theory of Hilbert
$C^\ast-$module and we denote it by $\ha$ (see \cite{lan}).
Obviously $\ha$ is countably generated and possesses an orthonormal
basis $\{e_i\otimes 1\}$, where $\{e_i\}$ is an orthonormal basis in
$H$. If $E$ is a countably generated Hilbert $\ma-$module then $E$
is unitarily equivalent to a fully complemented submodule of $\ha$
(see \cite{lan}). So we only consider Hilbert $C^\ast-$module $\ha$
in this paper.

We introduce a class of module maps which is analogous to rank-1
operators on a Hilbert space. $\forall x,y\in \ha$, define
$\theta_{x,y}:\ha\rightarrow\ha$ by $\theta_{x,y}(\xi)=\langle
\xi,y\rangle x$, $\forall \xi\in\ha$. Note that $\theta_{x,y}$ is
quite different from rank-1 linear operators on a Hilbert space to
some extents. For instance, we can not infer $x=0$ or $y=0$ from
$\theta_{x,y}=0$. But $\theta_{x,y}$ have the following properties:
$\theta_{x,\alpha y}=\theta_{\alpha^\ast x,y}$ and when $\ma$ is
commutative, $\alpha\theta_{x,y}=\theta_{\alpha x,y}$. We denote
$span_\ma\{\theta_{x,y}\}$ by $\mathcal{F}(\ha)$

In this paper we mainly consider module map
$\Phi:\mathcal{F}(\ha)\rightarrow \mathcal{F}(\ha)$ which maps
$\theta_{x,y}$ to some $\theta_{s,t}$. The method is analogous to
that of Hou's in \cite{hou1}, but much more complicated since many
properties in linear space can not generalized to the module
setting. Using our results we can calculate the free Fisher
information of a semicircular variable with rank one preserving
covariance.

\section{Rank one preserving module maps on $\mathcal{F}(\ha)$}

We first introduce a class of elements in $\ha$.

\begin{defn}
Let $\varepsilon$ be a orthonormal bases in $\ha$. $x\neq 0\in\ha$
will be called coordinatly invertible if $\langle e,x\rangle$ is
invertible unless $\langle e,x\rangle=0$, $\forall e\in\varepsilon$.
\end{defn}

Denote the set of all the coordinatly invertible elements in $\ha$
by $CI(\ha)$ or $CI$ for short. Obviously, $\varepsilon\subseteq
CI$.

We will see cordinatly invertible elements are similar to elements
in Hilbert space to some extents.

\begin{lem}
Let $y\in CI$ and $\theta_{x,y}=0$. Then $x=0$
\end{lem}
\begin{proof}
From $y\in CI$, there is $e\in\varepsilon $, such that $\langle
e,y\rangle\neq 0$ and invertible. Thus $\langle e,y\rangle x=0$ and
$x=0$.
\end{proof}

Note that the above lemma does not hold for general $x,y$.

The following lemma is well known in linear space and we can
generalize it to the modular setting.

\begin{lem}
Let $\mm$ be a Hilbert $\ma-$module, where $\ma$ is an unital
$C^\ast-$algebra, and let $\phi,\sigma:\mm\rightarrow \ma$ be
$\ma-$linear operators. Suppose $\sigma$ vanishes on the kernel of
$\phi$. Then there exists $b\in\ma$ such that $\sigma=\phi\cdot b$.
\end{lem}
\begin{proof}
Define $\widehat{\phi}: \mm/ker\sigma\rightarrow\ma$ by
$\widehat{\phi}(x+ker\phi)=\phi(x)$ and define
$\widehat{\sigma}:\mm/ker\sigma\rightarrow \ma$ by
$\widehat{\sigma}(x+ker \sigma)=\sigma(x)$. Then it is easy to see
$\widehat{\phi}$, $\widehat{\sigma}$ are $\ma-$linear and injective.

Now we write $x+ker\sigma=\widehat{\phi}^{-1}(a)$, where
$a=\widehat{\phi}(x+ker\sigma)=\phi(x)$. Then
$x+ker\sigma=a\widehat{\phi}^{-1}(1)$ and we have
\begin{eqnarray*}
&& \widehat{\sigma}(x+ker
\sigma)=\widehat{\sigma}(a\widehat{\phi}^{-1}(1))=a\widehat{\sigma}(\widehat{\phi}^{-1}(1))\\
&=&\widehat{\phi}(x+ker\sigma)\widehat{\sigma}(\widehat{\phi}^{-1}(1))
\end{eqnarray*}
So $\widehat{\sigma}=\widehat{\phi}\cdot
\widehat{\sigma}(\widehat{\phi}^{-1}(1))$ and $\sigma=\phi\cdot
\widehat{\sigma}(\widehat{\phi}^{-1}(1))$. We finish the proof by
letting $b=\widehat{\sigma}(\widehat{\phi}^{-1}(1))$.
\end{proof}

\begin{cor}
Let $g_1,g_2\in\mm$. If for all $x\in\mm$, $\langle x, g_1\rangle=0$
implying $\langle x, g_2\rangle=0$. Then there is $a\in\ma$, such
that $g_2=ag_1$.
\end{cor}
\begin{proof}
Define $\varphi_{g_i}(x)=\langle x, g_i\rangle,\, i=1,2$. Obviously
$\varphi_{g_i}s$ are $\ma-$linear. From Lemma 3, we know $\exists
b\in\ma$, such that $\varphi_{g_2}=\varphi_{g_1}\cdot b$, that is
$\langle x, g_2\rangle=\langle x,g_1\rangle b=\langle x, b^\ast
g_1\rangle$, $\forall x\in\mm$. Putting $a=b^\ast$, we get
$g_2=ag_1$.
\end{proof}

Now we consider $\ma-$linear operator $\Phi:
\mathcal{F}(\ha)\rightarrow \mathcal{F}(\ha)$, which satisfying
$\forall x\in\ha,y_0\in CI$, $\exists t_0\in CI$, such that
$\Phi(\theta_{x,y_0})=\theta_{s,t_0}$. Then we call $\Phi$ is Rank
decreasing. If $\forall x\neq 0$ implying $s\neq 0$, then $\Phi$
will be called Rank-1 preserving.

\begin{lem}
Let $\ma$ be a unital $C^\ast-$algebra. $x_1,x_2\in\mm,g_1,g_2\in
CI$ satisfying $\theta_{x_1,y_1}+\theta_{x_2,g_2}=\theta_{x_3,g_3}$.
Then at least one of the following cases occurs:\\
(i) $\exists\alpha_1\in\ma$, such that $g_1=\alpha_1 g_2$;\\
(ii) $\exists\alpha_2\in\ma$, such that $g_2=\alpha_2 g_1$;\\
(iii) $\exists \beta_1,\beta_2\in\ma$, such that $x_1=\beta_1 x_3$
and $x_2=\beta_2 x_3$.
\end{lem}

\begin{proof}
We will complete the proof by considering the following four cases.

Case (1) If $\forall \xi\in\ha$, $\langle\xi, g_1\rangle=0$ implying
$\langle\xi, g_2\rangle=0$. From Corollary 4, there exists
$\alpha_1\in\ma$ such that $g_1=\alpha_1 g_2$.

Case (2) If $\forall \xi\in\ha$, $\langle\xi, g_1\rangle=0$ implying
$\langle \xi,g_2\rangle=0$. Still from Corollary 4, there exists
$\alpha_2\in\ma$ such that $g_2=\alpha_2 g_1$.

Case (3) $\exists \xi_0\in\ha$, such that
$\langle\xi_0,g_2\rangle=0$ but $\langle \xi_0, g_1\rangle\neq 0$.
We can find $e\in\varepsilon$, such that $\langle e,g_2\rangle =0$
but $\langle e, g_1\rangle\neq 0$. Then we get
$$\langle e,g_1\rangle x_1+\langle e,g_2\rangle x_2=\langle
e,g_3\rangle x_3,$$ and
$$\langle e,g_1\rangle x_1=\langle e, g_3\rangle x_3$$
Since $g_1\in CI$, we have $x_1=\langle e,g_1\rangle^{-1}\langle e,
g_3\rangle x_3$. We put $\beta_1=\langle e, g_1\rangle ^{-1}\langle
e,g_3\rangle$. Then
$\theta_{\beta_1x_3,g_1}+\theta_{x_2,g_2}=\theta_{x_3,g_3}$ and thus
$\theta_{x_2,g_2}=\theta_{x_3,g_3-\beta_1^\ast g_1}$. Now we can
find $e'\in\varepsilon$ such that $x_2=\langle
e',g_2\rangle^{-1}\langle e',g_3-\beta_1^\ast g_1\rangle x_3$.
Letting $\beta_2=\langle e',g_2\rangle^{-1}\langle
e',g_3-\beta_1^\ast g_1\rangle$, then we get  (iii).

Case (4) $\exists \xi_0\in\ha$, such that
$\langle\xi_0,g_1\rangle=0$ but $\langle \xi_0, g_2\rangle\neq 0$.
Similar to Case (3), we get (iii) again.
\end{proof}

From the proof of Lemma 5, we have the following corollary.

\begin{cor}
With the notations in the above lemma, supposing $g_1\neq\alpha g_2$
and $g_2\neq \beta g_1$, $\forall \alpha,\beta\in\ma$ and $g_3\in
CI$, then $\beta_1$ or $\beta_2$ can be chosen to be invertible.
\end{cor}

We introduce some new notations. $\forall x,y\in\ha$,
$L^{CI}_x:=\{\theta_{x,g}\mid g\in CI\}$;
$R_y^{CI}:=\{\theta_{h,y}\mid h\in CI\}$; $L_x:=\{\theta_{x,g}\mid
g\in \ha\}$; $R_y:=\{\theta_{h,y}\mid h\in \ha\}$.

\begin{lem}
$\Phi$ is a rank decreasing  $\ma-$linear map. Then one of the
following cases holds:\\
(i) $\forall x\in\ha$, $\exists y\in\ha$, such that
$\Phi(L_x^{CI})\subseteq L_y^{CI}$;\\
(ii) $\forall x\in \ha$, $\exists f\in\ha$, such that
$\Phi(L_x^{CI})\subseteq R_f$
\end{lem}

\begin{proof}
If $\forall x\in\ha$, $\exists x_0,g_0$ such that
$\Phi(L_x^{CI})=\alpha \theta_{x_0,g_0}$, where $\alpha\in \ma$,
then (i), (ii) both hold. So we mainly consider the case which such
$x_0, g_0$ do not exist.

Assume there exists $x_{00}\in\ha$ such that
$\Phi(L^{CI}_{x_{00}})\nsubseteq L_x^{CI}$,
$\Phi(L_{x_{00}}^{CI})\nsubseteq R_f$, $\forall x,f\in\ha$. Then
there are $f_1,f_2\in CI$, such that
$\Phi(\theta_{x_{00},f_1})=\theta_{x_1,g_1}$ and
$\Phi(\theta_{x_{00},f_2})=\theta_{x_2,g_2}$ where $x_1\neq\alpha_1
x$, $x_2\neq \alpha_2 x, \forall x\in\ha, \alpha_1,\alpha_2,\in\ma$
and $g_1\neq \beta_1 g_2, g_2\neq\beta_2g_1,\forall
\beta_1,\beta_2\in\ma$.

Write $\Phi(\theta_{x_{00},f_1+f_2})=\theta_{x_3,g_3}$. From Lemma 5
$\theta_{x_1,g_1}+\theta_{x_2,g_2}\neq \theta_{x_3,g_3}$, $\forall
g_3,x_3\in\ha$ and this is a contradiction.

Thus we have proved $\forall x\in\ha$, either
$\Phi(L_x^{CI})\subseteq L_y^{CI}$ or $\Phi(L_x^{CI})\subseteq R_f$.

Now we take $x_0\in CI$ and suppose $\Phi(L_{x_0}^{CI})\subseteq
L_{y_0}^{CI}$. $M:=\{x\in\ha\mid \Phi(L_x^{CI})\subseteq
L_{y(x)}^{CI}, \hbox{ for some } y(x)\in\ha\}$ and $N:=\{x\in\ha\mid
\Phi(L_x^{CI})\subseteq R_{f(x)}, \hbox{ for some } f(x)\in CI
\hbox{ and }\Phi(\theta_{x,h})\neq \alpha \theta_{y, f(x)}\hbox{ for
all } h\in CI, y\in\ha\}$

Obviously, $M\bigcup N=\ha$, $M\bigcap N=\emptyset$.

Assuming $N\neq\emptyset$, then there exists $x_1\in N$,  such that
$\Phi(\theta_{x_1,f})=\theta_{y_1(f),g_1}$,for all $f\in CI$. On the
other hand $\Phi(\theta_{x_0,f})=\theta_{y_0,g_0(f)}$. Suppose
$x_0+x_1\in N$, then $\Phi(\theta_{x_0+x_1,f})=\theta_{y_2(f),g_2}$
and thus
\begin{equation}
\theta_{y_0,g_0(f)}+\theta_{y_1(f),g_1}=\theta_{y_2(f),g_2}
\end{equation}
Since $x_1\in N$, we can choose $f_0\in CI$,
$\Phi(\theta_{x_1,f_0})=\theta_{y_1(f_0),g_1}$ such that
$y_1(f_0)\neq \alpha y_0,\, \forall \alpha\in\ma$. From (1),
\begin{equation}
\theta_{y_0,g_0(f_0)}+\theta_{y_1(f_0),g_1}=\theta_{y_2(f_0),g_2}
\end{equation}

Assuming there exists $\alpha_0\in\ma$ such that $y_2(f_0)=\alpha_0
y_0$ and then
\begin{equation}
\theta_{y_0,g_0(f_0)}+\theta_{y_1(f_0),g_1}=\theta_{\alpha_0y_0,g_2}
\end{equation}
\begin{equation}
\theta_{y_1(f_0),g_1}=\theta_{y_0,\alpha_0^\ast g_2-g_0(f_0)}
\end{equation}
Since $g_1\in CI$, there is $e\in \varepsilon$ such that $\langle e,
g_1\rangle$ is invertible. Then we have
$$y_1(f_0)=\langle e, g_1\rangle ^{-1}\langle e, \alpha^\ast
g_2-g_0(f_0)\rangle y_0$$ which contradicting to $y_1(f_0)\neq
\alpha y_0$ and thus $y_2(f_0)\neq \beta y_0,\, \forall
\beta\in\ma$.

Then from lemma5, $g_2=\alpha_0 g_1$, and $\forall f\in CI$
\begin{equation}
\theta_{y_0,g_0(f)}+\theta_{y_1(f),g_1}=\theta_{y_2(f), \alpha_0g_1}
\end{equation}
Then we have
\begin{equation}
\theta_{y_0,g_0(f)}=\theta_{\alpha_0^\ast y_2(f)-y_1(f),g_1}
\end{equation}

For $e\in\varepsilon$,
\begin{equation}
\langle e, g_0(f)\rangle y_0=\langle e,g_1\rangle [\alpha^\ast_0
y_2(f)-y_1(f)]
\end{equation}

 Let $\beta_f y_0=\alpha_0^\ast y_2(f)-y_1(f)$ where $\beta_f\in\ma$
 and then from (6), we get
 \begin{equation}
 \theta_{y_0,g_0(f)}=\theta_{y_0,\beta_f^\ast g_1}.
 \end{equation}
 Since $x_0\in CI$ implying $y_0\in CI$, we have
 $g_0(f)=\beta_f^\ast g_1$ which contradicting to the assumption in
 the beginner of our proof. This is shown that $N=\emptyset$ i.e.
 $\ha=M$.

 If $\Phi(L_{x_0}^{CI})\subseteq R_f$, then by the similar
 method we can show $N=\ha$ and the proof is finished.
 \end{proof}

 \begin{prop}
 Let $\ma$ be a unital commutative $C^\ast-$algebra and
  let $\Phi$ be a rank decreasing $\ma-$module map. Then one of
  the following cases occurs\\
  i) $\forall x\in\ha$, $\Phi(L_x)\subseteq L_y$\\
  ii) $\forall x\in\ha$, $\Phi(L_x)\subseteq R_f$.
  \end{prop}

  \begin{proof}
  $\forall x,g\in\ha$, then $g=\sum\limits_i \alpha_i g_i$, where
  $\alpha_i\in\ma$ and $g_i\in CI$. Suppose $\Phi(L_x^{CI})\subseteq
  L_y^{CI}$ and we have
  \begin{eqnarray*}
  &&\Phi(\theta_{x,g})=\Phi(\theta_{x,\sum\limits_i \alpha_i
  g_i})=\sum\limits_i\alpha_i^\ast \Phi(\theta_{x,g_i})\\
  &=&\sum\limits_i \alpha_i^\ast \theta_{y,h_i}=\theta_{y,\sum\limits_i
  \alpha_i h_i}
  \end{eqnarray*}
  Thus we infer that $\Phi(L_x)\subseteq L_y$.
  \end{proof}

  In order to characterize $\Phi$, we also need the following lemma
  which is well known in linear space.

  \begin{lem}
  Let $\ma$ be a unital commutative $C^\ast-$algebra and let $A$ be
  an injective $\ma-$linear map on $\ha$. There are $x_1,x_2\in \ha$
  such that $ Ax_1\neq \alpha Ax_2$ and $Ax_2\neq \beta Ax_1$, $\forall
  \alpha,\beta\in\ma$. $B$ is another $\ma-$linear map satisfying
  that $\forall x\in\ha$, there exists $\lambda_x$ such that $Bx=\lambda_x
  Ax$. Then $B=\lambda A$, for some $\lambda\in\ma$.
  \end{lem}
  \begin{proof}
  We will complete the proof by two steps.

  (1) Suppose that there are $x_1,x_2\in\ha$, such that $\langle
  Ax_1,Ax_2\rangle=0$ (In fact, we can choose $x_1\in CI$).

  From the condition, $\exists \lambda_1,\lambda_2,\lambda_3\in\ma$,
  such that $Bx_1=\lambda_1Ax_1$, $Bx_2=\lambda_2 Ax_2$, and $B(x_1+x_2)=\lambda_3
  A(x_1+x_2)$ and so
  \begin{equation*}
  (\lambda_1-\lambda_3)Ax_1+(\lambda_2-\lambda_3)Ax_2=0.
  \end{equation*}
Then
\begin{eqnarray*}
&&\langle(\lambda_1-\lambda_3)Ax_1+(\lambda_2-\lambda_3)Ax_2,
(\lambda_1-\lambda_3)Ax_1+(\lambda_2-\lambda_3)Ax_2\rangle\\
&=&(\lambda_1-\lambda_3)\langle Ax_1,Ax_1\rangle
(\lambda_1-\lambda_3)^\ast +(\lambda_2-\lambda_3)\langle
Ax_2,Ax_2\rangle (\lambda_2-\lambda_3)^\ast\\
&=& 0.
\end{eqnarray*}
Thus
\begin{equation*}
(\lambda_1-\lambda_3)\langle Ax_1,Ax_1\rangle
(\lambda_1-\lambda_3)^\ast=0;
\end{equation*}
\begin{equation*}
(\lambda_2-\lambda_3)\langle Ax_2,Ax_2\rangle
(\lambda_2-\lambda_3)^\ast=0;
\end{equation*}
\begin{equation*}
(\lambda_1-\lambda_3)Ax_1=(\lambda_2-\lambda_3)Ax_2=0.
\end{equation*}
From the above equations, we get
\begin{equation*}
Bx_1=\lambda_1 Ax_1=\lambda_3 Ax_1;
\end{equation*}
\begin{equation*}
Bx_2=\lambda_2Ax_2=\lambda_3Ax_2.
\end{equation*}

$\forall x\in \ha$ such that $Ax,Ax_1,Ax_2$ is a orthogonal set
($x_1\in CI)$. We claim that $\exists \lambda\in\ma$ such that
$Bx=\lambda Ax,\,\forall x\in\ha$. In fact
\begin{equation*}
Bx=\lambda_x Ax=\lambda_{x+x_1} Ax
\end{equation*}
\begin{equation*}
Bx_1=\lambda_1 Ax_1=\lambda_{x+x_1}Ax_1=\lambda_3 Ax_1
\end{equation*}
Since $A$ is injective, we get $\lambda_{x+x_1}x_1=\lambda_3x_1$,
i.e. $(\lambda_{x+x_1}-\lambda_3)x_1=0$. $\forall e\in\varepsilon$,
such that $\langle x_1,e\rangle\neq 0$, and since $x_1\in CI$, we
have $(\lambda_{x+x_1}-\lambda_2)\langle x_1,e\rangle=0$ and
$\lambda_{x+x_1}=\lambda_3$. So $Bx=\lambda_3 Ax$.

(2) Let $Ax_1,Ax_2,\cdots$ be the orthonormal bases in $A(\ha)$.
$\forall x\in\ha$, $Ax=\alpha_1 Ax_1+\alpha_2 Ax_2+\cdots$ and since
$A$ is injective, we have $x=\alpha_1x_1+\alpha_2 x_2+\cdots$. Then
\begin{eqnarray*}
Bx &=& B(\alpha_1 x_1+\alpha_2 x_2+\cdots)\\
&=& \alpha_1 B(x_1)+\alpha_2 B(x_2)+\cdots\\
&=&\alpha_1\lambda Ax_1+\alpha_2\lambda Ax_2+\cdots\\
&=&\lambda(\alpha_1 Ax_1+\alpha_2 Ax_2+\cdots)\\
&=&\lambda Ax
\end{eqnarray*}
\end{proof}

\begin{prop}
Let $\ma$ be an unital commutative $C^\ast-$algebra and let
$\Phi:\ha\rightarrow \ha $ be a rank-1 preserving $\ma-$linear.
$\forall x,f\in\ha$, $\Phi(L_x^{CI})\subseteq L_{\varphi(x)}^{CI}$,
$\Phi(R_f^{CI})\subseteq R_{r(f)}^{CI}$. There are $f_1,f_2\in CI$
such that $r(f_1)\neq \alpha r(f_2)$, $r(f_2)\neq \beta r(f_1)$,
$\forall \alpha,\beta\in\ma$. Then there exist $\ma-$linear maps
$A,C:\ha\rightarrow \ha$ such that
$\Phi(\theta_{x,f})=\theta_{A(x),C(f)}$, $\forall x,f\in\ha$.
\end{prop}

\begin{proof}
It follows from Proposition 8 that $\forall f\in CI$, there is a map
$\psi_f$ on $\ha$ such that
\begin{equation}
\Phi(\theta_{x,f})=\theta_{\psi_f(x),r(f)}
\end{equation}
We will complete the proof by 5 steps.

{\bf Step 1.} We show $\psi_f$ is injective. If not, there are
$x_1\neq x_2\in\ha$, such that $\psi_f(x_1)=\psi_f(x_2)$. Then
\begin{equation*}
\Phi(\theta_{x_1,f})=\theta_{\psi_f(x_1),r(f)}
\end{equation*}
\begin{equation*}
\Phi(\theta_{x_1,f})=\theta_{\psi_f(x_2),r(f)}
\end{equation*}
and thus $\Phi(\theta_{x_1-x_2,f})=\theta_{0,r(f)}=0$ which
contradicting to $\Phi$ preserving "rank-1".

{\bf Step 2.} We show $\psi_f$ is $\ma-$linear. $\forall x,y\in\ha$,
\begin{equation*}
\Phi(\theta_{x+y,f})=\theta_{\psi_f(x+y),r(f)}=\theta_{\psi_f(x)+\psi_f(y),r(f)}
\end{equation*}
from lemma 1, we get $\psi_f(x+y)=\psi_f(x)+\psi_f(y)$. On the other
hand, since $\ma$ is commutative, $\theta_{\alpha
x,f}=\alpha\theta_{x,f}$ and therefore $\psi_f(\alpha
x)=\alpha\psi_f(x)$.

{\bf Step 3.} We show $\forall f\in CI$, $\exists$ $\ma-$linear map
$\psi$ such that $\psi_f=\alpha(f)\psi$, where $\alpha(f)\in\ma$.

It is easy to see we can choose $f_1,f_2\in CI$ satisfying
$r(f_1)\neq\alpha r(f_2),r(f_2)\neq r(f_1)$, $\forall
\alpha,\beta\in\ma$ and $f_1+f_2$ still in $CI$.

We consider
\begin{equation}
\Phi(\theta_{x,f_1+f_2})=\theta_{\psi_{f_1+f_2}(x),r(f_1+f_2)}=\theta_{\psi_{f_1}(x),r(f_1)}
+\theta_{\psi_{f_2}(x),r(f_2)}.
\end{equation}
Since $r(f_1)\neq \alpha r(f_2),r(f_2)\neq \beta r(f_1)$, from Lemma
5, $\exists \alpha_x, \beta_x\in\ma $ such that
\begin{equation*}
\psi_{f_1}(x)=\alpha_x\psi_{f_1+f_2}(x)
\end{equation*}
\begin{equation*}
\psi_{f_2}(x)=\beta_x\psi_{f_1+f_2}(x)
\end{equation*}
From Corollary 6 $\alpha_x$, or $\beta_x$ can be chosen to be
invertible.

Then follows from Lemma 9 and its proof we know there are
$\alpha_0,\beta_0\in\ma$ such that
\begin{equation*}
\psi_{f_1}=\alpha_0 \psi_{f_1+f_2}
 \end{equation*}
 \begin{equation*}
\psi_{f_2}=\beta_0 \psi_{f_1+f_2}
\end{equation*}
where $\alpha_0$ or $\beta_0$ can be invertible.

{\bf Claim.} $\forall f\in CI$, either $r(f)\neq \alpha r(f_1),
r(f_1)\neq \beta r(f)$ or $r(f)\neq \alpha r(f_2), r(f_2)\neq \beta
r(f)$.

In fact assume $\exists \alpha_{00},\beta_{00}\in\ma$ such that
$r(f)=\alpha_{00} r(f_1)$ or $r(f_1)=\beta_{00} r(f)$ we can show
$r(f)\neq \alpha r(f_2),r(f_2)\neq \beta r(f)$. If not
$r(f)=\alpha_1 r(f_2)$ or $r(f_2)=\beta_1 r(f)$. Supposing
$r(f)=\alpha_0 r(f_1)$ and $r(f)=\alpha_1 r(f_2)$ hold at the same
time then $\alpha_0$ is invertible since $r(f),r(f_1),r(f_2)\in CI$
and $r(f_1)=\alpha_0^{-1} \alpha_1 r(f_2)$ which contradicting to
the condition. So the claim has been shown.

Now for all $f\in CI$, suppose $r(f)\neq \alpha r(f_1),r(f_1)\neq
\beta r(f)$, there are $\alpha',\beta'$ such that
\begin{equation}
\psi_{f_1}=\alpha'\psi_{f+f_1}=\alpha_0\psi_{f_1+f_2}
\end{equation}
$\forall x_0\in CI$, we get
\begin{equation*}
\psi_{f_1}(x_0)=\alpha'
\psi_{f+f_1}(x_0)=\alpha_0\psi_{f_1+f_2}(x_0)
\end{equation*}
and we infer from
$\psi_{f_1+f_2}(x_0),\psi_{f+f_1}(x_0),\psi_{f_1}(x_0)\in CI$ that
$\alpha'$ is invertible. On the other hand
\begin{equation*}
\psi_f=\beta'\psi_{f+f_1}=\beta' \alpha'^{-1}\alpha_0\psi_{f_1+f_2}
\end{equation*}
Letting $\psi=\psi_{f_1+f_2}$ the desired result is obtained.

{\bf Step 4.} We show $\Phi(\theta_{x,f})=\theta_{A(x),C(f)}$,
$\forall f\in CI$, where $A,C$ are $\ma-$linear.

Let $A=\psi$ and $C_0(f)=\alpha(f) r(f)$. Then $\forall f\in CI$,
\begin{equation*}
\Phi(\theta_{x,f})=\theta_{A(x),C_0(f)}
\end{equation*}

{\bf Step 5.} For a general $f\in\ha$, $f=\sum\limits_i \beta_i
e_i$, where $e_i\in \varepsilon\subseteq CI$, we get
\begin{eqnarray*}
&&\Phi(\theta_{x,f})=\Phi(\theta_{x,\sum\limits_i \beta_i^\ast e_i})\\
&=&\sum\limits_i \beta_i \theta_{A(x),C_0(e_i)}=\sum\limits_i
\theta_{A(x),\sum\limits_i \beta^\ast_i C_0(e_i)}\\
&=&\theta_{A(x),C(f)}
\end{eqnarray*}
where $C(f):=\sum\limits_i \beta^\ast_i C_0(e_i)$.

When we take $x\in CI$,  it is easy to show $C$ is $\ma-$linear.
\end{proof}

By the similar way, we also have
\begin{cor}
If $\forall x,f\in\ha$, $\Phi(L_x^{CI})\subseteq R_{f(x)}$,
$\Phi(R_f^{CI})\subseteq L_{r(f)}$, then there are $\ma-$conjugate
linear maps $A,C$ such that $\Phi(\theta_{x,f})=\theta_{Af,Cx}$.
\end{cor}

Our main result in this section is
\begin{thm}
  $\Phi$ preserve rank 1 if and only if  $\Phi$ has one of the following  forms

  (1) $\exists$ injective $\ma-$ linear maps $A,C$ such that $\Phi(\theta_{x,f})=\theta_{Ax,Cf}$

  (2) $\exists$ injective $\ma-$ conjugate linear maps $A,C$ such
  that $\Phi(\theta_{x,f})=\theta_{Af,Cx}$.
\end{thm}

\begin{proof}
(1) $\forall x\in\ha$, $\Phi(L_x)\subseteq L_{\varphi(x)}$. One can
 easily infer that  $\Phi(R_f)\subseteq R_{f(x)},\forall f\in\ha$. Since $\Phi$ preserving
rank 1, there are $f_1,f_2\in CI$ such that $r(f_1)\neq \alpha
r(f_2), r(f_2)\neq \beta r(f_1)$. Then from Proposition 10 there
exist $\ma-$linear maps $A,C$ such that
$\Phi(\theta_{x,f})=\theta_{Ax,Cf}$. It is easy to observe that
$A,C$ are injective.

(2) $\forall x\in \ha$, $\Phi(L_x)\subseteq R_{f(x)}$. We also
obtain that $\Phi(R_f)\subseteq L_{r(f)}$. There are
$x_1,x_2\in\ma$, such that $\varphi(x_1)\neq \alpha \varphi(x_2)$,
$\varphi(x_2)\neq\beta\varphi(x_1)$. Then from Lemma again there
exist $\ma-$ conjugate linear maps $A,C$ such that
$\Phi(\theta_{x,f})=\theta_{Af,Cx}$.
\end{proof}

\section{Rank preserving module maps on $L(\ha)$}

The most important class of operators on Hilbert $C^\ast-$module is
adjointable operators. $T$, a operator on a Hilbert $C^\ast-$module,
will be called adjointable if there is a operator $S$ such that
$\langle Tx,y\rangle=\langle x, Sy\rangle$. $S$ is often denoted by
$T^\ast$. The set of all the adjointable operators on $\mathcal{H}$
will denoted by $L(\mathcal{H})$.

Let $\mathcal{H}$ be a Hilbert $\ma-$module. the strict topology on
$\mathcal{H}$ is defined by the family of seminorms $v\mapsto
\|\langle x, v\rangle \|, x\in\mathcal{H}$. From \cite{lan}, we know
that $\mathcal{F}(\ha)$ is strictly dense in $L(\ha)$.
\begin{thm}
Let $\Phi:L(\ha)\rightarrow L(\ha)$ be strictly continuous and
preserve rank 1 which has the form
$\Phi(\theta_{x,y}=\theta_{Ax,Cy}$ or
$\Phi(\theta_{x,y})=\theta_{Ay,Cx}$ If $A,C\in L(\ma)$ the $\Phi$
has one of the following forms:

(1) There are $A,B\in L(\ha)$ which are injective such that $\forall
T\in L(\ha)$, $\Phi (T)=ATB$;

(2) There are adjointable conjugate $\ma-$linear operators $A,B$
such that $\forall T\in L(\ha)$, $\Phi(T)=AT^\ast B$.

\end{thm}

\begin{proof}
Since $\Phi$ preserving rank-1, there exist injective module maps
$A,C$ such that $\Phi(\theta_{x,y})=\theta_{Ax,Cy}$ or there are
conjugate $\ma-$linear maps $A,C$ such that
$\Phi(\theta_{x,y})=\theta_{Ay,Cx}$.

We only consider the first case. $\Phi
(\theta_{x,y})=\theta_{Ax,Cy}=A\theta_{x,y}B$ where $B:=C^\ast$.
From \cite{lan}, we know that $\mathcal{F}(\ha)$ is dense in
$L(\ha)$ in the sense of strict topology. We have
\begin{equation*}
\Phi(T)=\Phi(\sum\limits_{i=1}^\infty
\theta_{x_i,y_i})=\sum\limits_{i=1}^\infty
\Phi(\theta_{x_i,y_i})=\sum\limits_{i=1}^\infty A\theta_{x_i,y_i}
B=ATB.
\end{equation*}

The other case can be proved similarly.
\end{proof}

\begin{cor}
With the notations and assumptions as the above theorem, if $\Phi$
is surjective then $A,B$ are invertible.
\end{cor}

\begin{proof}
Since $A,C$ are injective, $B$ is surjective. It follows from $\Phi$
surjective that $A,B$ are invertible.
\end{proof}

Free probability is a noncommutative probability theory. This
theory, due to D.Voiculescu, has very important applications on
operator algebras. In this section we mainly consider
operator-valued free probability theory. Let $\mm$ be a unital
algebra and $\mb$ be a subalgebra of $\mm$, $1\in \mb$, and let
$E:\mm\rightarrow \mb$ be a conditional expectation. We call
$(\mm,E,\mb)$ an operator-valued (or $\mb-$ valued) noncommutative
probability space and elements in $\mm$ are called random variables.
We can use cumulant function to describe a random variable (see
\cite{spe}) The most important class of random variables in free
probability is the semicircular variables. In an operator-valued
noncommutative probability space a semicircular variable is
connected with a linear map.
\begin{defn}\cite{spe}\cite{voi}
Let $(\mm,E,\mb)$ be a noncommutative probability space and let
$\eta:\mb\rightarrow\mb$ be a linear map. A self-adjoint element
$X\in\mm$ will be called a semicircular variable with covariance
$\eta$ (or $\eta-$semicircular variable) if it satisfies
$k^{(1)}(X)=0$, $k^{(2)}(X\otimes bX)=\eta(b)$, $k^{(m+1)}(X\otimes
b_1X\otimes\cdots\otimes b_mX)=0$, for all
$b,b_1,\cdots,b_m\in\mb,m\geq 2$, where $(k^{(n)})_{n\geq 1}$ is the
cumulant function induced by $E$.
\end{defn}

In \cite{meng2}, we have generalized the notion of free Fisher
information to the operator-valued setting.
\begin{defn}
Let $(\mm,E,\mb)$ be a $\mb-$valued noncommutative probability
space, and let $X\in\mm$ be a self-adjoint random variable.
$\eta:\mb\rightarrow \mb$ be a linear map. $\xi\in L^2(M)$ will be
called the conjugate variable of $X$ with respect to $\eta$, if it
satisfies: $k^{(1)}(\xi)=0$, $k^{(2)}(\xi\otimes bX)=\eta(b)$,
$k^{(m+1)}(\xi\otimes b_1X\otimes\cdots\otimes b_mX)=0$, $\forall
b,b_1,\cdots,b_m\in\mb,m\geq 2$.

Let $\tau$ be a faithful state on $\mb$.   The free Fisher
information of $X$ is defined by $\varphi^\ast_\tau
(X:\mb,\eta)=\tau E(\xi\xi^\ast)$.
\end{defn}

We usually use $J(X:L(\ha),\eta)$ to denote the conjugate variable
of $X$ with respect to $\eta$. To construct the conjugate variable
of a random variable is not easy in general. In
\cite{meng1},\cite{meng2}, we have calculated the free Fisher
information of a semicircular variable with conditional expectation
covariance. Now we calculate the free Fisher information of a
semicircular variable with rank-1 preserving covariance.

\begin{thm}
Let $\Phi$ be a surjective, rank-1 preserving $\ma-$linear map.
$\Phi(T)=ATB$, where $A,B\in L(\ha)$. $\eta:T\rightarrow B^{-1}TB$.
$X\in (\mm,E,L(\ha))$ be $\Phi-$semicircular. Then
$J(X:L(\ha),\eta)=XA^{-1}B^{-1}$. Let $\tau$ be a faithful tracial
state on $L(\ha)$, then
$\varphi_\tau^\ast(X:L(\ha),\eta)=\tau(B^{-1\ast}A^{-1\ast})$.
\end{thm}

\begin{proof}
We only need to verify $XA^{-1}B^{-1}$ satisfying the formulae in
Definition 15. Since $X$ is $\Phi-$semicircular, we have
$k^{(1)}(XA^{-1}B^{-1})=0$;
\begin{eqnarray*}
&&k^{(2)}(XA^{-1}B^{-1}\otimes bX)=E(XA^{-1}B^{-1}bX)\\
&=&AA^{-1}B^{-1}bB=B^{-1}bB=\eta(b);
\end{eqnarray*}
\begin{equation*}
k^{(m+1)}(XA^{-1}B^{-1}\otimes b_1X\otimes\cdots \otimes b_mX)=0
\end{equation*}
for all $b,b_1,\cdots,b_m\in\mb$, $m\geq 2$.

Thus $J(X:L(\ha),\eta)=XA^{-1}B^{-1}$.
\begin{eqnarray*}
\varphi^\ast_\tau(X:L(\ha),\eta)&=&\tau
E(J(X:L(\ha),\eta)J(X:L(\ha),\eta)^\ast)\\
&=&\tau E(XA^{-1}B^{-1}B^{-1\ast}A^{-1\ast} X)\\
&=& \tau\Phi (A^{-1}B^{-1}B^{-1\ast}A^{-1\ast} B)\\
&=&\tau (B^{-1\ast} A^{-1\ast})
\end{eqnarray*}
\end{proof}

\bibliographystyle{amsplain}

\begin{thebibliography}{{voi11}}

\bibitem{hou2}G.An, J.Hou, Rank-preserving multiplicative maps on
B(X), Lin. Alg. Appl. 342 (2002) 59-78

\bibitem{hou1}J.Hou, Rank preserving linear maps on B(X), Sci. in
China (ser. A), 32 (1989), 929-940

\bibitem{kap}I.Kaplansky, Modules over operator algebras, Amer. J.
Math. 75 (1953) 839-853

\bibitem{lan} C.Lance, Hilbert $C^\ast-$modules: a toolkit for
operator algebraists, London Math. Soc. Lecture Notes Series, vol.
210, Cambridge University Press, Cambridge, 1994

\bibitem{meng1}B.Meng, M.Guo, X.Cao, Operator-valued free Fisher
information and modular frame, Proc. Amer. Math. Soc., vol 133,
no. 10, 3087-3096

\bibitem{meng2}B.Meng, M.Guo, X.Cao, Some applications of free
Fisher information on frame theory, J. Math. Anal. Appl. 311
(2005) 466-478

\bibitem{spe} R.Speicher, Combinatorial theory of the free product
with amalgamation and operator-valued free probability theory, Mem.
Amer. Math. Soc., Vol. 627, 1998

\bibitem{voi} D.Voiculescu, Operations on certain non-commutative
operator-valued random variables, Ast\'{e}risque 32 (1995) 243-275













\end{thebibliography}

\end{document}